\documentclass{article}
\usepackage{graphicx}
\usepackage{amsmath}
\usepackage{amsfonts}
\usepackage{amssymb}
\usepackage{amsthm}
\usepackage{mathtools}
\usepackage{enumerate}
\usepackage{lscape}
\usepackage{longtable}
\usepackage{rotating}
\usepackage{multirow}
\usepackage{url}
\usepackage{subfigure}
\usepackage{rotating}
\usepackage{comment}
\usepackage[utf8]{inputenc}
\usepackage{amssymb}
\usepackage{amsfonts}
\usepackage{tikz}
\usepackage{natbib} 
\usepackage{wrapfig}
\usepackage{color}
\usepackage{hyperref}

\usepackage{algorithm}
\usepackage{algpseudocode}

\theoremstyle{definition}

\title{A Three-Tier Time-Scale Architecture for Controlling Complex Nonlinear Systems}
\author{Vyacheslav Kungurtsev}
\date{December 2025}

\begin{document}
\maketitle

\begin{abstract}
This letter proposes a three-tier computational architecture for the real-time control of nonlinear complex systems, such as time-dependent PDEs. There is an important class of such problems for which existing single- and two-time-scale approaches are fundamentally insufficient due to lack of a priori system knowledge, computational complexity, model fidelity requirements, and uncertainty. The proposed architecture consists of an offline, meso-scale, and real-time layer of computation, with distinct roles for each layer and specific information flow between them. The result is a practical systems-level paradigm that enables real-time operation of complex nonlinear control problems.
\end{abstract}

\begin{sloppypar}

\section{Introduction}

We propose a three-tier computational architecture organized by time scale and model fidelity that enables real-time control under severe computational and uncertainty constraints. The intention is to provide a generalized framework for the real-time operation of very complex nonlinear systems, such as evolution Partial Differential Equations (PDEs) and nonlinear (e.g. chaotic) dynamics. The properties of such state dynamics present serious issues of intractability of operation with sufficient efficiency and performance when addressed solely at real-time computation, or even with two-time scale schemes. 

In a full-length companion draft,~\cite{difonzo2025towards}, a specific concrete operating problem is presented. Specifically, a nonlinear hyperbolic partial differential equation modeling open channel water flaw - the shallow water equations in particular, are taken as state-defining constraints in a Dynamic Programming problem. Motivated by control for hydropower in service of energy demand, real-time operation with awareness and robustness with respect to future uncertainty is considered. This problem is a desiderata of water engineering that could significantly improve hydropower operations. However, the significant fundamental, methodological, and practical challenges present require a comprehensive computational systems engineering architecture to reasonably address this problem.

In the next Section, the computational systemic reasoning for three tiers at three distinct time scales is presented. Subsequently, the generalized formal framework of computation, latency, and information exchanged between levels is stated in Section~\ref{sec:formal}. The modular incorporation of existing architectures and tools into the three-tier framework is presented in Section~\ref{sec:relate}. Concluding remarks comprise Section~\ref{sec:conc}.

This paper intends to provide the motivation for and a useful set of guidelines for a systems engineering approach towards computational operation of complex nonlinear systems. The companion manuscript~\cite{difonzo2025towards}, which fleshes out a comprehensive grounding of the time scale procedure, illustrates the generic systems engineering principles presented in this letter. As a first-principles contribution to computational systems engineering, the paper seeks to provide a usable and effective framework to approach challenging real-world optimization tasks. 

\textbf{Contributions}:
\begin{enumerate}
    \item We formalize the necessity of three computational time scales for real-time control of complex nonlinear systems.
    \item We define the role, latency, and information flow of offline, meso, and real-time layers.
    \item We relate the architecture to existing control, optimization, data assimilation, and scientific computing frameworks.
\end{enumerate}

\section{Why Three Tiers?}\label{sec:why}
One of the primary challenges with many nonlinear PDEs or other dynamical systems is the lack of a priori knowledge of the control-to-state map and function space regularity. As a result, the range of available computational procedures becomes delicate and sparse, while being insufficient by themselves to solve the problem. Fast solution techniques that are simultaneously reliable becomes intractable to any known software or algorithms, \emph{by themselves}. However, high fidelity (e.g.) CFD simulations are available for the systems, nonlinear programming solvers are available to solve discretized optimal controls, and real-time control procedures operate when the system is sufficiently simple and low-dimensional. By an agglomeration of computational effort across these toolkits, adequately integrated, a sensible approach can be developed for ambitious control tasks.

The theoretical challenges of concern are discussed for the specific case of optimal control with the Shallow Water Equations for hydropower operation in the companion draft~\cite{difonzo2025towards}. To begin with, the nonlinear hyperbolic PDE does not admit any global existence theory and no uniqueness guarantees. This implies that any control to state map may or may not exist for an arbitrary value of the control, and exhibit poor sensitivity with respect to disturbances in the problem data. Second, the space of functions, potentially a Banach Algebra or other general space, may not be amenable to adjoint computation or numerical methods with approximation guarantees. Finally, the overall nonlinearity, mixed-integer, and uncertainty-awareness complications of the problem, especially with any desiderata to achieve close to global optimality, is computationally intractable with respect to complexity scaling considerations. 

At real-time, i.e., seconds to minutes of latency, only small and simple optimization problems can be solved. ``Real-Time Iteration'' (RTI)~\cite{diehl2005real} solves quadratic programs meant to approximate an NLP representing an evolving Model Predictive Control (MPC) problem. The incorporation of strongly nonlinear dynamics and fine spatial discretization needed for PDE approximation make real-time Newton-based approaches not suitable. At the real-time level quadratic and linear optimization problems are solved, with mixed integer considerations only recently been explored~\cite{quirynen2021sequential}. 

Consider then a two-tier structure. At a slower level of computational latency, optimization problems with a complete and precise discretization of the PDE or other evolution/dynamic system, as well as the uncertainty, is meant to be solved. Model Order Reduction (MOR) is used to communicate solutions at this slower computational operation to the real-time layer. Real-time updates are communicated to the slower level, that are then used within a parametric path-following approach to incorporate within the nonlinear programming (NLP) based software design. However, with the complex problems of interest, even this two level structure would not be sufficient. Large scale optimization of nonlinear PDEs typically requires many hours if not days of computation on an HPC cluster. Scenario generation multiplies this required load. Moreover, when it is unknown that the equation defining the control's influence on the state exhibits a solution, and  one with smoothness properties amenable to computation, the reliable operation of NLP is far from guaranteed at reasonable latency. 

Thus, we need an additional third level of computational scaffolding. The highest level is meant to be performed \emph{offline}, with essentially unlimited time and any available use of computational HPC resources. The aim of the offline layer is to perform a significant quantity of computation that prepares a two-level system for temporal operation. The offline layer computes a comprehensive catalog of solution estimates at various conditions as well as an empirical understanding of the control-to-state map availability. This catalog is then used by the \emph{meso} level of the two-level temporal operation to prepare candidate solutions for its numerical optimization. Regions of decision and state space that have been found to exhibit solutions with reasonable numerical sensitivity offline present an operational surface on which optimization at the meso layer can be performed. Information flowing from the real-time to the meso to the offline level, as far as the outcomes of various operating scenarios, completes the full three-tier loop with feedback. This defines the comprehensive end-to-end Dual-Control (real-time estimation and control) computational systems architecture.

\section{Formal Time-Scale Separation and Information Flow}\label{sec:formal}

\paragraph{The Offline Layer} is meant to use maximal HPC resources available and the time necessary to perform its primary role in the system. In practice, this can be from a year up to three years of dedicated computation while the rest of the layers are undergoing research and development. The generous computation presents the opportunity to traverse the problem's complexity explosion, as far as PDE space-time discretization, uncertainty simulation, and the search for (approximate) global optimality. The \emph{Primary Goal} of the offline layer is to create 1) a means of quickly identifying a region of the decision space wherein a control-to-state map locally exists and is regular and 2) a catalog of solutions that can be used as reference trajectories and warm starts for embedded-in-time operation. 

\emph{Communication} By frontloading extensive computation before the rollout of a time-embedded control scheme, the layer delivers a useful database and associated query tools that the meso layer can utilize for its operation. Namely, the database stores probabilistic information on the well-behavedness of the control-to-state as well as a large collection of reference solutions. The meso layer can query the database in order to steer optimization as well as provide candidate solutions and warm starts.

\paragraph{The Meso Layer}
The meso layer is time-embedded, meaning there is some latency $T_2$ that is targeted as far as the duration of an iteration in its control algorithm. This duration is not fast enough for desired real-time operation, but works to provide a more accurate and precise background computation that steadily communicates to the real-time layer. In practice $T_2$ can be a few minutes to a few hours. The meso layer performs space-time discretizations of the underlying Dynamic Programming problem, and applies off-the-shelf NLP solvers to compute solutions. A set of candidate solution ``particles'' can be maintained. Some are sequentially dropped when poor performance or constraint violation is predicted, replaced by off-line catalog solutions or random mutation of other particles. The best candidate solution(s) will be communicated to the real-time layer. 

\emph{Communication} The Meso layer periodically communicates updated optimal solution(s) to the real-time layer. These define the optimal control sequence when incorporating the nonlinear dynamics comprehensively in the model used for optimization as well as the uncertainty.

The Meso layer queries the Offline layer's catalog in order to restrict the decision space to compact sets wherein a control-to-state map is likely to exist and be smooth as well as to provide candidate reference solutions. Meanwhile, information on the discrepancy between predicted and actual behavior, informed by real-time sensor data passed from the real-time layer, is then passed on to the offline layer. This information can be used in the offline layer to adaptively improve its catalog and decision regularity classification map.

\paragraph{The Real-Time Layer} is, finally, meant to operate at the latency $T_1<T_2$ desired by the operator for online operation. In practice, this can be a few seconds up to minutes. In general, one or two orders of magnitude, in seconds, is meant to separate the real-time and meso layer. In real-time only small to medium-scale quadratic programs, possibly mixed-integer, are solved. Real-time iteration algorithms use Newton type methods to trace a nonlinear solution. This layer will perform the operative MPC scheme - it iteratively solves an optimization problem with updated information and implements the control that is found to be the solution for the current time. Fast latency presents circumstances where by design and operation the optimization problems don't significantly change across each time step, facilitating fast computation through warm and hot starts.

\emph{Communication} The real-time layer periodically receives information regarding new nonlinear solutions from the meso layer. In order to facilitate integration, MOR is performed in the meso layer to communicate accurate but computationally amenable reference solutions. The real-time layer can use these solutions as tracking objectives, presenting convexity in the objective functional. The real-time layer, meanwhile, communicates information obtained by real-time operating sensors with streaming data on the state and the environment. This information is sent to the meso layer, which also performs data mining to send refined information to the offline layer. 

\begin{wrapfigure}{r}{0.5\textwidth}
\begin{center}\includegraphics[scale=0.4]{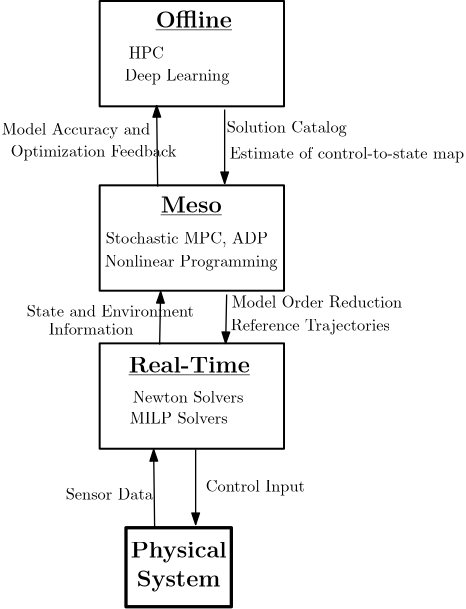}
    \caption{Summary of the three time scale computational system architecture.}\label{fig}
    \end{center}
\end{wrapfigure}
\paragraph{Summary} Figure~\ref{fig} presents the architecture. This operation - offline catalog to meso NLP solutions to real-time fast computation, and back up with real-time information, describes a comprehensive closed loop engineering system.


\section{Relation to Existing Architectures }\label{sec:relate}
In this Section, for each layer, as well as the cross-layer communication, a suggested toolkit is provided for the application of existing algorithms and software to facilitate the operation. This serves to motivate the understanding of the purpose and role of the computational system components and provide a reference for practitioners. In order to avoid a long reference list, citations will not be included. However, most of these methods are mentioned in the companion manuscript~\cite{difonzo2025towards}, which presents a comprehensive grounding of the 3-tier computational system architecture for the case study of PDE control for water engineering. 

\textbf{Offline} Mixed Integer Nonlinear Programming, Multistage Stochastic Programming, Reinforcement Learning, Physics Informed Machine Learning, Distributed and shared memory highly parallel high performance computing

\textbf{Meso} Stochastic Model Predictive Control, Approximate Dynamic Programming, Scenario Generation, Nonlinear Programming, Particle Evolution Strategies

\textbf{Real-Time} Real-time Iteration, Advanced-step MPC, Newton Solvers, Quadratic Program Solvers, Mixed Integer Linear/Quadratic Program Solvers, Embedded Control

\textbf{Meso to Real-time Integration} Model Order Reduction, Fine to Coarse grid interpolation and MGOpt,

\textbf{Real-time to Meso Feedback} Moving Horizon Estimation, Filtering Methods, Embedded Sensing, System Identification

\textbf{Offline to Meso Integration} Initial guesses, Catalog for mutation selection in particle evolution, Bayesian methods (for uncertainty quantification of the control-to-state map)

\textbf{Meso to Offline Feedback} Data mining, Reinforcement Learning based Model Predictive Control, Uncertainty Quantification

\section{Conclusion}\label{sec:conc}
Real and potential aspirations of real-time control of complex nonlinear systems present a formidable challenge for computational engineering. Even a two time scale operating system would be inadequate for this problem for many cases of interest. This paper presents a three-time-scale computational engineering system, with specific roles for each layer and integration protocols, intended to serve as an operating paradigm for such problems, and is hopefully a useful reference for researchers, engineers and other stakeholders in complex systems control design. 
\paragraph*{Acknowledgements}
This work received funding from the National Centre for Energy II (TN02000025). The author would like to thank Fabio DiFonzo, Michael Holst, Sven Leyffer, Paul Manns, Wladimir Neves and Thomas Surowiec for valuable insight and discussions on the problem.

\end{sloppypar}

\bibliographystyle{plain}
\bibliography{refs}

\begin{thebibliography}{1}

\bibitem{diehl2005real}
Moritz Diehl, Hans~Georg Bock, and Johannes~P Schl{\"o}der.
\newblock A real-time iteration scheme for nonlinear optimization in optimal
  feedback control.
\newblock {\em SIAM Journal on control and optimization}, 43(5):1714--1736,
  2005.

\bibitem{difonzo2025towards}
Fabio DiFonzo, Michael Holst, Morteza Kimiaei, Vyacheslav Kungurtsev, and
  Songqiang Qiu.
\newblock Towards real time control of water engineering with nonlinear
  hyperbolic partial differential equations.
\newblock {\em arXiv preprint arXiv:2512.14387}, 2025.

\bibitem{quirynen2021sequential}
Rien Quirynen and Stefano Di~Cairano.
\newblock Sequential quadratic programming algorithm for real-time
  mixed-integer nonlinear mpc.
\newblock In {\em 2021 60th IEEE Conference on Decision and Control (CDC)},
  pages 993--999. IEEE, 2021.

\end{thebibliography}

\end{document}